\def\Arxived{arxived}
\newif\ifarxived
\ifx\Arxived\Kaba\else\arxivedtrue\fi
\newif\ifextended
\ifx\Extended\Kaba\else\extendedtrue\fi
\newif\ifprivate
\ifx\Private\Kaba\else\privatetrue\extendedtrue\fi
\newif\ifJapanese
\newif\iftesting
\ifextended\testingtrue\fi
\ifarxived
\extendedtrue
\testingfalse
\privatefalse
\fi
\documentclass[12pt]{article}
\ifextended
\ifarxived
\usepackage[bookmarks=true,bookmarksnumbered=true,
  bookmarkstype=toc,
  pdftitle={Generically supecompact cardinals by forcing with chain conditions},%
  pdfauthor={Sakaé Fuchino (\quad\quad\quad), Hiroshi Sakai (\quad\quad\quad)},%
  pdfsubject={},%
  pdfkeywords={generically supercompact cardinal, ccc, weakly Mahloness, tree property}]{hyperref}
\else
\usepackage[dvipdfmx,bookmarks=true,bookmarksnumbered=true,
  bookmarkstype=toc,
  pdftitle={Generically supecompact cardinals by forcing with chain conditions},%
  pdfauthor={Sakaé Fuchino (\qquad\qquad\quad), Hiroshi Sakai (\qquad\qquad\quad)},%
  pdfsubject={},%
  pdfkeywords={generically supercompact cardinal, ccc, weakly Mahloness, tree property}]{hyperref}
\usepackage[dvipdfmx]{pxjahyper}
\hypersetup{
  colorlinks=true,
  linkcolor=blue,
  filecolor=blue,      
  citecolor=cyan,
  urlcolor=cyan
}
\fi 
\else 
\newcommand{\href}[2]{#2}
\def\hyperref[#1]#2{#2}
\newcommand{\phantomsection}{}
\fi

\ifarxived
\usepackage{graphicx, color}
\else
\usepackage[dvipdfmx]{graphicx, color}
\fi
\usepackage{amsmath, amssymb}
\usepackage{bbm}
\usepackage{dsfont}
\usepackage{bbold}
\usepackage{marginnote}
\usepackage{mathtools}
\usepackage[normalem]{ulem}
\usepackage{setspace}

\definecolor{darkelectricblue}{rgb}{0.33, 0.39, 0.52}
\definecolor{darkgreen}{rgb}{0.31, 0.47, 0.26}
\newcommand{\extendedcolor}{\color{darkelectricblue}}
\newcommand{\darkgreen}{\color{darkgreen}}
\newcommand{\privatecolor}{\color{darkgreen}}
\ifextended
\newcommand{\darkred}{\color[rgb]{0.8,0.1,0.1}}
\newcommand{\It}{\it\darkred{}}
\else
\newcommand{\darkred}{}
\newcommand{\It}{\it{}}
\fi
\usepackage{theorem}
\newcommand{\bbd}[1]{{\mathbb{#1}}}
\theorembodyfont{\ifJapanese\rm\else\it\fi}
\newcount\minute	
\newcount\hour		
\newcount\hourMins  
\ifJapanese
\def\today%
{
  \the\year /\,\zeroPadTwo{\the\month}/\,\zeroPadTwo{\the\day}%
}
\fi
\def\now%
{
  \minute=\time    
  \hour=\time \divide \hour by 60 
  \hourMins=\hour \multiply\hourMins by 60
  \advance\minute by -\hourMins 
  \zeroPadTwo{\the\hour}:\zeroPadTwo{\the\minute}%
}
\def\zeroPadTwo#1%
{
  \ifnum #1<10 0\fi    
  #1
}
\title{\scalebox{0.94}[1.6]{\bf Generically supercompact cardinals}\\[2\jot]
\scalebox{0.94}[1.6]{\bf by forcing with chain conditions}}
\author{\protect\scalebox{1}[1.4]{\ \qquad\qquad\quad\qquad\qquad\qquad\quad \qquad\qquad\qquad\qquad}\medskip\\
  \protect\scalebox{1}[1.4]{\quad Saka\'e Fuchino$^{*,\dagger}$\qquad\quad\ Hiroshi Sakai$^{*,\ddagger}$}\\
  }
\date{}

\renewcommand{\baselinestretch}{1.2}
\renewcommand{\thefootnote}{(\arabic{footnote})\,}
\iftesting
\fi
\iftesting
\newcommand{\Label}[1]{\label{#1}\marginnote{{\renewcommand{\baselinestretch}{0.4}\tiny 
		  #1}}}
\newcommand{\LabelR}[1]{\label{#1}\marginnote{{\renewcommand{\baselinestretch}{0.4}\tiny 
		  \hfill\rlap{#1}}}}
\else
\newcommand{\Label}[1]{\label{#1}}
\newcommand{\LabelR}[1]{\label{#1}}
\fi

\def\memo#1{\ifprivate\marginnote{{\darkgreen\normalsize%
      \renewcommand{\baselinestretch}{0.4}\tiny\mbox{}\vspace{-2.52ex}
      \par\relax%
			#1\par\mbox{}}}\else\fi}%
\newcounter{frml}[section]
\newcounter{frmla}[section]
\def\thefrml{{\arabic{section}.\arabic{frml}}}
\def\thefrmla{{$\aleph$\arabic{section}.\arabic{frmla}}}
\def\frmlabel#1{\refstepcounter{frml}{\def\baka{#1}\ifx\baka\empty\else\label{#1}\fi}%
{\rm({\thefrml})\hfill\hfill\hfill}}
\def\frmlabela#1{\refstepcounter{frmla}{\def\baka{#1}\ifx\baka\empty\else\label{#1}\fi}%
{\rm({\thefrmla})\hfill\hfill\hfill}}
\def\xitem[#1]{\item[\frmlabel{#1}]\mbox{}%
	\iftesting\marginnote{{\renewcommand{%
				\baselinestretch}{0.6}\tiny#1}}\fi\ignorespaces}
\def\xitemq[#1]{\item[\frmlabel{#1}]\mbox{}%
	\ignorespaces}
\def\xitemd[#1]#2{\item[(\ref{#1})$#2$\hfill\hfill\hfill]}
\def\xitemA[#1]{\item[\frmlabela{#1}]\mbox{}%
	\iftesting\marginnote{{\renewcommand{%
				\baselinestretch}{0.6}\tiny#1}}\fi\ignorespaces}
\def\xitemsub[#1]#2{\item[\frmlabel{#1}$_{#2}$]\mbox{}%
	\iftesting\marginnote{{\renewcommand{%
				\baselinestretch}{0.6}\tiny#1}}\fi\ignorespaces}
\def\xxitem[#1][#2]{\item[(\ref{#1}{\makebox[1.4ex][c]{#2}})]\mbox{}%
	\iftesting\marginnote{{\renewcommand{%
				\baselinestretch}{0.6}\tiny\{#1\}\{#2\}}}\fi\ignorespaces}
\def\xitemof#1{{\rm({\ref{#1}})}}
\def\Xitem[#1]{\item[{\makebox[7ex][l]{\rm(\ref{#1})}}]\iftesting\marginnote{{\renewcommand{%
				\baselinestretch}{0.6}\tiny#1}}\fi\ignorespaces}

\newenvironment{xitemize}{\begin{list}{}{\parsep=0.5\smallskipamount%
			\itemindent=-0.4ex%
			\itemsep=0.5\smallskipamount\leftmargin=4em\labelwidth=3em\labelsep=0.7em}}%
							 {\end{list}}
\def\assert#1{\noindent\makebox[4.8ex][r]{\rm(\makebox[2.2ex][c]{#1})}\ \ \ignorespaces}
\def\wassert#1{\assert{#1}}
\def\wassertof#1{\makebox[4.8ex][r]{\rm(\makebox[2.2ex][c]{#1})}}%
\def\assertof#1{(#1)}%

\def\daimaru#1{\makebox[1em][c]{\mbox{\leavevmode\lower.144ex\hbox{
        \rlap{\hbox to 
          0.76em{\hfil\mbox{}\hfill{}\raisebox{0.054ex}{\scalebox{1.2}{○}}\hfil}}
        \raise0.342ex\hbox to 1em{\hfil{\hspace{0.16em}\footnotesize#1}\hfil}}}}\,}

\newtheorem{Thm}{\ifJapanese{\bf 定理}\else {\bf Theorem}\fi}[section]

\newtheorem{ThmA}{\ifJapanese{\bf 定理\,A\!}\else{\bf Theorem\,A\!}\fi}[section]
\newtheorem{Ex}[Thm]{\ifJapanese{\bf 例}\else {\bf Example}\fi}
\newtheorem{Prop}[Thm]{\ifJapanese{\bf 命題}\else{\bf Proposition}\fi}

\newtheorem{Lemma}[Thm]{\ifJapanese{\bf 補題}\else{\bf Lemma}\fi}
\newtheorem{LemmaA}[ThmA]{\ifJapanese{\bf 補題\,A\!}\else{\bf Lemma\,A\!}\fi}
\newtheorem{Cor}[Thm]{\ifJapanese{\bf 系}\else{\bf Corollary}\fi}

\newtheorem{Claim}{{\bf Claim}}[Thm]

\newcommand{\prf}{\ifJapanese{\bf 証明．\ }\ignorespaces\else{\bf 
		Proof.\ \ }\ignorespaces\fi}

\newcommand{\prfofClaim}{\raisebox{-.4ex}{\Large $\vdash$\ \ }}

\newcommand{\prfof}[1]{\ifJapanese{\bf #1 の証明．\ \ }%
	\ignorespaces\else{\bf Proof of #1:}\ \ \ignorespaces\fi}
\newcommand{\Thmof}[1]{\ifJapanese{定理\,\ref{#1}}\else{Theorem~\ref{#1}}\fi}

\newcommand{\Lemmaof}[1]{\ifJapanese{補題\,\ref{#1}}\else{Lemma~\ref{#1}}\fi}

\newcommand{\LemmaAof}[1]{\ifJapanese{補題\,A\,\ref{#1}}\else{Lemma\,A\,\ref{#1}}\fi}

\newcommand{\Propof}[1]{\ifJapanese{命題\,\ref{#1}}\else{Proposition~\ref{#1}}\fi}

\newcommand{\Claimof}[1]{{Claim \ref{#1}}}

\newcommand{\Corabove}{{\ifJapanese 系\else Corollary\fi\ \number\theThm}}

\newcommand{\Claimabove}{{Claim \number\theClaim}}

\newcommand{\ubecause}[3]{\underbrace{{}#1{}%
  \ifx\bakakaba#2\bakakaba\rule[-0.72ex]{0pt}{1pt}\else\rule[#2]{0pt}{1pt}\fi}_{\footnotesize\clap{#3}}}
\newcommand{\obecause}[3]{\overbrace{{}#1{}%
  \ifx\bakakaba#2\bakakaba\rule[1.62ex]{0pt}{1pt}\else\rule[#2]{0pt}{1pt}\fi}^{\footnotesize\clap{#3}}}
\newsavebox{\qedbox}\sbox{\qedbox}{
{\unitlength=0.05mm \begin{picture}(40,60)
\put(0,0){\framebox(30,44)[cc]{}}
\put(30,-7){\rule{7\unitlength}{44\unitlength}}
\put(10,-7){\rule{27\unitlength}{7\unitlength}}
\end{picture}}}
\newcommand{\qed}{\mbox{}\hfill\usebox{\qedbox}}
\newcommand{\smallqed}%
{\mbox{}\smallskip\hfill\raisebox{-.4ex}{\Large $\dashv$}}
\newcommand{\qedof}[1]%
{\mbox{} \hspace*{\fill}{\usebox{\qedbox}{\tiny~(#1)}}}
\newcommand{\Qedof}[1]%
{\mbox{} \hspace*{\fill}{\usebox{\qedbox}%
{\tiny~(#1~\number\theThm)}}}
\newcommand{\QedAof}[1]%
{\mbox{} \hspace*{\fill}{\usebox{\qedbox}%
{\tiny~(#1~\number\theThmA)}}}
\newcommand{\qedofThm}{\Qedof{\ifJapanese 定理\else Theorem\fi}}

\newcommand{\qedofCor}{\Qedof{\ifJapanese 系\else Corollary\fi}}
\newcommand{\qedofProp}{\Qedof{\ifJapanese 命題\else Proposition\fi}}
\newcommand{\qedofLemma}{\Qedof{\ifJapanese 補題\else Lemma\fi}}
\newcommand{\qedofLemmaA}{\QedAof{\ifJapanese 補題\,A\!\else Lemma\,A\!\fi}}

\newcommand{\qedskip}{\medskip}

\newcommand{\qedofClaim}%
{\mbox{}\hfill\raisebox{-.4ex}{\Large $\dashv$ }\nolinebreak%
\mbox{\tiny~(Claim~\number\theClaim)}}
\newcommand{\qedofClaimA}%
{\mbox{}\hfill\raisebox{-.4ex}{\Large $\dashv$ }\nolinebreak%
\mbox{\tiny~(Claim~A\,\number\theClaimA)}}
\newcommand{\qedofClaimAof}[1]%
{\mbox{}\hfill\raisebox{-.4ex}{\Large $\dashv$ }\nolinebreak%
\mbox{\tiny~(Claim~A\,\ref{#1})}}
\newcommand{\qedofSubclaim}%
{\mbox{}\hfill\raisebox{-.4ex}{\Large $\dashv$ }\nolinebreak%
\mbox{\tiny~(Subclaim~\number\theSubclaim)}}
\newcommand{\cardof}[1]{\mathopen{|\,}#1\mathclose{\,|}}
\newcommand{\ccardof}[1]{\mathopen{\,\|}#1\mathclose{\|\,}}
\newcommand{\ulsetof}[1]{\mathopen{\,|}#1\mathclose{|\,}}
\newcommand{\Card}{{\it Card\/}}
\newcommand{\Reg}{{\it Reg\/}}
\newcommand{\setof}[2]{\{#1\,:\,#2\}}
\newcommand{\ssetof}[1]{\{#1\}}

\newcommand{\subseteqand}[1]{\mathrel{\mathop{\subseteq}%
		\limits_{\scriptscriptstyle\hbox to 14pt{$\scriptscriptstyle #1$\hss}}}}

\newcommand{\mapping}[3]{#1:#2\rightarrow #3}

\newcommand{\Elembed}[4]{#1:#2\stackrel{\prec\hspace{0.8ex}}{\rightarrow}_{#4}#3}

\newcommand{\fnsp}[2]{\mbox{}^{{#1}\hspace{-0.02em}}#2}
\newcommand{\imageof}{{}^{\,{\prime}{\prime}}}

\newcommand{\seqof}[2]{\langle#1\,:\,#2\rangle}
\newcommand{\pairof}[1]{\langle#1\rangle}
\newcommand{\psof}[1]{{\mathcal P}\/(#1)}
\newcommand{\forces}[2]{\,\|\hspace{-.35ex}\mbox{\sf--}_{\,#1\,}%
\mbox{\rm``}\,#2\,\mbox{\rm''}}

\newcommand{\modelof}[1]{\models\!\mbox{\rm``\,}#1\mbox{\rm''}}

\newcommand{\crit}{\mbox{\it crit\/}}

\newcommand{\circleq}{\mathrel{{\leqslant}%
		\hspace{-0.86ex}{\lower-0.53ex\hbox{$\scriptscriptstyle\circ$}}}}

\newcommand{\restr}{\restriction}
\newcommand{\id}{\mathop{id\,}}

\newcommand{\cf}{\mathop{cf\/}}

\newcommand{\mahlo}{\mathrm{M\ell\,}}

\newcommand{\sk}{\mathop{\it sk}}
\newcommand{\Fn}{{\rm Fn}}

\newcommand{\dom}{\mathop{\it dom}}

\newcommand{\poP}{\bbd{P}}
\newcommand{\poQ}{\bbd{Q}}

\newcommand{\On}{{\rm On}}

\newcommand{\genG}{\mathbb{G}}

\newcommand{\genH}{\mathbb{H}}

\newcommand{\condp}{\mathbbm{p}}

\newcommand{\condr}{\mathbbm{r}}

\newcommand{\LT}{{<}\,}

\newcommand{\GT}{{>}\,}
\newcommand{\GE}{{\geq}\,}

\newcommand{\ctenten}{,\mbox{}\hspace{0.08ex}{.}{.}{.}\hspace{0.1ex}}
\newcommand{\ctentenc}{,{}\linebreak[0]\hspace{0.04ex}{{.}{.}{.}\hspace{0.1ex},\,}\linebreak[0]}

\newcommand{\xmbox}[1]{ $\relax{\rm #1}\relax$ }
\newcommand{\strA}{\mathfrak{A}}
\newcommand{\strB}{\mathfrak{B}}
\newcommand{\strC}{\mathfrak{C}}

\newcommand{\continuum}{2^{\aleph_0}}

\newcommand{\calC}{{\mathcal C}}
\newcommand{\calD}{{\mathcal D}}

\newcommand{\calF}{{\mathcal F}}
\newcommand{\calG}{{\mathcal G}}
\newcommand{\calH}{{\mathcal H}}
\newcommand{\calI}{{\mathcal I}}
\newcommand{\calJ}{{\mathcal J}}

\newcommand{\calM}{{\mathcal M}}

\newcommand{\calP}{{\mathcal P}}

\newcommand{\calS}{{\mathcal S}}

\newcommand{\varin}{\mathrel{\varepsilon}}

\newcommand{\ZFC}{{\sf ZFC}}

\newcommand{\MA}{{\sf MA}}

\newcommand{\FRP}{{\sf FRP}}

\newcommand{\RC}{{\sf RC}}

\newcommand{\refl}{{\mathfrak{r}\mathfrak{e}\mathfrak{f}\mathfrak{l}\,}}

\newcommand{\REFL}{{\mathfrak{R}\mathfrak{E}\mathfrak{F}\mathfrak{L}\,}}

\newcommand{\st}{such that}
\newcommand{\wrt}{with respect to}
\newcommand{\Wolog}{Without loss of generality}
\newcommand{\wolog}{without loss of generality}

\newcommand{\tfae}{the following are equivalent}

\newcommand{\po}{poset}
\newcommand{\pos}{posets}
\newcommand{\uniV}{{\sf V}}

\newcommand{\Pkl}[2]{\ifx\bakakaba#1\bakakaba\ifx\bakakaba#2\bakakaba{\mathcal 
    P}_\kappa(\lambda)\else{\mathcal P}_\kappa(#2)\fi\else{\mathcal P}_{#1}(#2)\fi}

\newcommand{\utildeT}[1]{%
  \hbox to 0pt{\smash{$\mathop{\textstyle #1}\limits_{%
			\raisebox{0.4ex}[0pt]{$\scriptstyle\sim$}}$}\hss}%
  \relax\phantom{\mathord{{#1}_{\rule[-0.6ex]{0pt}{1pt}}}}}
\newcommand{\utildeS}[1]{%
	\hbox to 0pt{\smash{$\mathop{\scriptstyle #1}\limits_{%
				\raisebox{0.6ex}[0pt]{$\scriptscriptstyle\sim$}}$}\hss}%
	\relax\phantom{\mathord{{#1}_{\rule[-0.6ex]{0pt}{1pt}}}}}
\newcommand{\utildeSS}[1]{%
	\hbox to 0pt{$\mathop{\scriptscriptstyle #1}%
		\limits_{\scriptscriptstyle\sim}$\hss}%
		\relax\phantom{\underline{#1}}}
\newcommand{\utilde}[1]{%
	\mathchoice{\utildeT{#1}}{\utildeT{#1}}{\utildeS{#1}}{\utildeSS{#1}}}

{\end{minipage}\end{trivlist}}

\begin{document}
\maketitle
\renewcommand{\thefootnote}{$\ast$\ }
  \footnotetext{Graduate School of System Informatics, Kobe University \\Rokko-dai 1-1, Nada, Kobe 657-8501 Japan
   \\
    \quad\scalebox{0.95}[1]{\tt ${}^\dagger$ fuchino@diamond.kobe-u.ac.jp, ${}^\ddagger$ hsakai@people.kobe-u.ac.jp}}

\ifextended
\phantomsection
\addcontentsline{toc}{section}{* Generically supercompact cardinals by forcing with chain conditions}
\addcontentsline{toc}{section}{**** by S.Fuchino and H.Sakai}
\fi

\ifextended
\addcontentsline{toc}{section}{Abstract}
\fi
\begin{abstract}
A ccc-generically supercompact cardinal $\kappa$ can be smaller than or equal to the continuum. On the other 
hand, such a cardinal $\kappa$ still satisfies diverse largeness properties, like that it is a 
stationary limit of ccc-generically measurable cardinals (\Thmof{P-gen-5}). This is in a strong 
contrast to $\calP$-generically supercompact cardinals for the class $\calP$ of 
all $\sigma$-closed \pos, which can 
be $\aleph_n$ for any $n>1$. 
\end{abstract}

\ifextended
{\extendedcolor 
\addcontentsline{toc}{section}{Contents}
\newcommand{\myscalebox}[1]{\scalebox{0.88}[1.06]{#1}}
\begin{quotation}
	\footnotesize
	\noindent
	\centerline{
      \normalsize\tt\quad\ Contents\hspace{6em}\mbox{}}\mbox{}\\
     {\mbox{}\hspace{-1.6em}\tt\makebox[3.4ex][l]{\ref{intro}.}%
      \hyperref[intro]{\tt\myscalebox{Introduction and preliminaries}}}\ \ \dotfill\ \ {\pageref{intro}}\\
     {\mbox{}\hspace{-1.6em}\tt\makebox[3.4ex][l]{\ref{large}.}%
      \hyperref[large]{\tt\myscalebox{Generically weakly compact cardinals}}}%
      \ \ \dotfill\ \ {\pageref{large}}\\        
     {\mbox{}\hspace{-1.6em}\tt\makebox[3.4ex][l]{\ref{meas}.}%
      \hyperref[meas]{\tt\myscalebox{Generically measurable cardinals}}}%
      \ \ \dotfill\ \ {\pageref{meas}}\\        
     {\mbox{}\hspace{-1.6em}\tt\makebox[3.4ex][l]{\ref{gen-sc}.}%
      \hyperref[gen-sc]{\tt\myscalebox{Generically supercompact cardinals}}}%
      \ \ \dotfill\ \ {\pageref{gen-sc}}\\        
     {\mbox{}\hspace{-1.6em}\tt\makebox[3.4ex][l]{\ref{refl}.}%
      \hyperref[refl]{\tt\myscalebox{Reflection properties down to $\LT$ a generically supercompact cardinal}}}%
      \ \ \dotfill\ \ {\pageref{refl}}\\

      \noindent
     {\mbox{}\hspace{-1.6em}\hyperref[ref]{\tt
      References}}\ \ \dotfill\ \ {\pageref{ref}}\\ 

\end{quotation}}
\fi
\renewcommand{\thefootnote}{}
\footnotetext{{\it Date:} May 9, 2021
  \qquad {\it Last update:} 
  \today\ (\now\ JST)\vspace{-1\smallskipamount}
}
\footnotetext{{\it MSC2020 Mathematical Subject Classification:}
  03E35, 03E50, 03E55, 03E57, 03E65\vspace{-1\smallskipamount}}
\footnotetext{\mbox{}\\[-2.52ex]The research is supported by Kakenhi Grant-in-Aid for Scientific Research (C) 20K03717}

\ifextended
\ifprivate
\footnotetext{This is a private extended version of a draft of the paper with the same title to appear in 
  the RIMS K\^oky\^uroku volume on RIMS set theory workshop 2021.  
  All additional 
  details not to be contained the submitted version of the paper are either typeset in 
  {\extendedcolor "dark electric blue"} in case the text is also included in the extended 
  version of the paper or in {\darkgreen "dark green"} if the comment should appear only in 
  the private version. 
  The numbering of the assertions is kept identical with the submitted version.

  The most up-to-date file of this private extended version is downloadable as:\\
  \href{https://fuchino.ddo.jp/papers/RIMS2021-ccc-gen-supercompact-xx.pdf}{\tt 
  https://fuchino.ddo.jp/papers/RIMS2021-ccc-gen-supercompact-xx.pdf} 
}
\else
\footnotetext{\extendedcolor This is an extended version of a draft of the paper with the same title to appear in 
  the RIMS K\^oky\^uroku volume on RIMS set theory workshop 2021.  
  All additional 
  details not to be contained \extendedcolor in the submitted version of the paper are either typeset in 
  dark electric blue (the color in which this paragraph is typeset) or put in separate appendices. 
  The numbering of the assertions is kept identical with the submitted version.

  The most up-to-date file of this extended version is downloadable as:\\
  \href{https://fuchino.ddo.jp/papers/RIMS2021-ccc-gen-supercompact-x.pdf}{\tt 
  https://fuchino.ddo.jp/papers/RIMS2021-ccc-gen-supercompact-x.pdf} 
}
\fi
\else
\footnotetext{An updated and extended version of this paper with more details and 
  proofs is downloadable as:\qquad \href{https://fuchino.ddo.jp/papers/RIMS2021-ccc-gen-supercompact-x.pdf}{\tt https://fuchino.ddo.jp/papers/RIMS2021-ccc-gen-supercompact-x.pdf}}
\fi

\renewcommand{\thefootnote}{\arabic{footnote})\,}
\section{Introduction and preliminaries.}
\Label{intro}
For a  class $\calP$ of \pos, we say that a cardinal $\kappa$ is {\It $\calP$-generically 
measurable}\/ ({\It$\calP$-g.\ measurable}, for short) if there is $\poP\in\calP$ \st, for a
$(\uniV,\poP)$-generic $\genG$, there  
are $j$, $M\subseteq\uniV[\genG]$ \st\ $\uniV[\genG]\models\Elembed{j}{\uniV}{M}{\kappa}$ 
holds\footnote{When we write $\Elembed{j}{N}{M}{\kappa}$, we mean
$N$ is a transitive model (possibly a class model) of $\ZFC^-$, 
$j$ is an elementary embedding of $\uniV$ into $M$, $M$ is transitive, and $\kappa$ is 
the critical point of $j$. }. If $\kappa$ is 
$\ssetof{\poP}$-generically measurable, we shall also say that $\kappa$ 
is $\poP$-generically measurable.

A cardinal $\kappa$ is {\It $\calP$-generically $\lambda$-supercompact}\/ ({\It$\calP$-g.\
  $\lambda$-supercompact}, for short) for a given cardinal $\lambda\geq\kappa$, 
for short) if 
there is $\poP\in\calP$ with a $(\uniV,\poP)$-generic $\genG$ and $j$,
$M\subseteq\uniV[\genG]$ with
\begin{xitemize}
\xitem[x-intro-0] $\uniV[\genG]\models\Elembed{j}{\uniV}{M}{\kappa}$, $j(\kappa)>\lambda$,
$j\imageof\lambda\in M$.
\end{xitemize}

A cardinal $\kappa$ is {\It $\calP$-generically supercompact}\/ ({\It$\calP$-g.\
  supercompact}, for short) if it is $\calP$-g.\
  $\lambda$-supercompact for all $\lambda\geq\kappa$.

Clearly, for $\kappa<\lambda<\lambda'$, $\calP$-g.\ $\lambda'$-supercompactness of $\kappa$ 
implies $\calP$-g.\ $\lambda$-super\-compact\-ness of $\kappa$ 
and $\calP$-g.\ $\kappa$-supercompactness of $\kappa$ is equivalent to $\calP$-g.\ 
measurability of $\kappa$. 

In the following we mainly consider the cases in which $\calP$ is the class of all $\nu$-cc \pos\ 
for some uncountable $\nu$. In 
this case, we shall say {\It $\nu$-cc-generically measurable} ({\It $\nu$-cc-g.\ measurable}, for short), 
or {\It $\nu$-cc-generically $\lambda$-supercompact}\/ ({\It $\nu$-cc-g.\ $\lambda$-supercompact}, for short), 
in place of $\calP$-generically measurable or $\calP$-generically $\lambda$-supercompact, respectively.

Starting from a measurable (supercompact, resp.) cardinal $\kappa$, it is easy to obtain a model 
where a ccc-g.\ measurable (supercompact, resp.) cardinal is less than or equal to the 
continuum. Actually, forcing with $\Fn(\lambda,2)$ for any $\lambda\geq\kappa$ will create 
such a model. 

We can also consider the generic versions of weak compactness: A cardinal $\kappa$ is said 
to be {\It $\calP$-generically weakly compact}\/ for a class $\calP$ of \pos\ (or {\It
  $\calP$-g.\ weakly compact}, for short), if, for any 
$A\subseteq\kappa$ ($A\in\uniV$), there is a transitive set model $M$ of $\ZFC^-$ with $\kappa$, $A\in M$ 
\st, for some $\poP\in\calP$ and $(\uniV,\poP)$-generic $\genG$, we have 
$\Elembed{j}{M}{N}{\kappa}$ for some $j$, $N\in\uniV[\genG]$. 

We shall also say {\It $\nu$-cc-g.\ weakly compact}\/ etc.\ similarly to above.

\begin{Lemma}
\Label{P-intro-0} For a class $\calP$ of \pos, if $\kappa$ is $\calP$-generically 
measurable then $\kappa$ is $\calP$-generically weakly compact. 
\end{Lemma}
\prf Suppose that $\poP\in\calP$ and $(\uniV,\poP)$-generic $\genG$ are \st\ there 
are $j^*$, $M^*\subseteq\uniV[\genG]$ with $\Elembed{j^*}{\uniV}{M^*}{\kappa}$. Let
$M:=\calH(\kappa^+)$, $N:=\calH(j^*({\kappa)^+}^{M^*})^{M^*}$ and 
$j:=j^*\restr M$. Then, these $M$, $\poP$, $\genG$, $j$, $N$ are witnesses of the property in 
the definition of the $\calP$-g.\ weakly compactness for all $A\subseteq\kappa$. 
\qedofLemma
\qedskip

We refer mainly \cite{millennium-book} for results in connection with precipitousness and generic 
ultrapower while our notation tend to be more compatible with that of \cite{higher-inf}. 
Names in forcing are denoted by alphabets with undertilde adopting the notation of 
\cite{proper-improper}.

\section{Generically weakly compact cardinals}
\Label{large}

\begin{Lemma}
\LabelR{P-ccc-gen-0} Suppose that $\kappa$ is $\nu$-cc-g.\ weakly compact for a $\nu<\kappa$. Then \smallskip

\wassert{1} $\kappa$ is weakly Mahlo.

\wassert{2} $\kappa$ has the tree property.
\end{Lemma}
\prf Assume that $\kappa$ is $\nu$-cc-g.\ weakly compact. 

\assertof{1}: \assertof{a}\quad First, we prove that $\kappa$ is not a successor 
cardinal\footnote{Actually, we can skip \assertof{a} since \assertof{b} implies 
  \assertof{c} and this establishes   
  \assertof{1} (see \Lemmaof{P-gen-0-0}). }. Suppose, toward a 
contradiction, that $\kappa$ is a successor cardinal, say $\kappa=\mu^+$. Note that
$\nu\leq\mu$. 

Let 
$A\subseteq\kappa$ be a set which codes $\seqof{s_\xi}{\xi<\kappa}$ where each $s_\xi$ for 
$0<\xi<\kappa$ is a surjection from $\mu$ to $\xi$.

Let $M$ be a transitive model of $\ZFC^-$ \st\ $\kappa$, $A\in M$ and there is a $\nu$-cc \po\
$\poP$ with $(\uniV,\poP)$-generic $\genG$ \st\ there are $j$, $N\in\uniV[\genG]$ \st\
\begin{xitemize}
\xitem[x-ccc-gen-0]  
  $\uniV[\genG]\models\Elembed{j}{M}{N}{\kappa}$.
\end{xitemize}

Now, since $A\in M$, we have $M\models\kappa=\mu^+$ and, since $j(\mu)=\mu$ 
by $\mu<\kappa$, we have $N\modelof{j(\kappa)=\mu^+}$ by elementarity. Thus
\begin{xitemize}
\item[]  $j(\kappa)=(\mu^+)^N\leq(\mu^+)^{V[\genG]}\obecause{=}{}{\qquad\qquad\qquad
  \vbox{\hbox{since $\poP$ preserves 
  cardinals $\GT\nu$ \vspace{-1\smallskipamount}}\hbox{by the $\nu$-cc of $\poP$}}}(\mu^+)^V=\kappa$. 
\end{xitemize}
This is a contradiction to $\kappa=\crit(j)$. 

\assertof{b}\quad 
Next, we prove that $\kappa$ is regular. Suppose, again toward a contradiction, that 
$\kappa$ is singular and let $\seqof{\kappa_\xi}{\xi<\delta}$ be a strictly increasing 
sequence of cardinals $<\kappa$ cofinal in $\kappa$ and \st\ $\delta<\kappa$. Let 
$A\subseteq\kappa$ be a set which codes the sequence $\seqof{\kappa_\xi}{\xi<\delta}$.

Let $M$ be a transitive model of $\ZFC^-$ \st\ $\kappa$, $A\in M$ and there is a $\nu$-cc \po\
$\poP$ with $(\uniV,\poP)$-generic $\genG$ \st\ there are $j$, $N\in\uniV[\genG]$ \st\
\xitemof{x-ccc-gen-0} holds.

By $A\in M$, we have $\seqof{\kappa_\xi}{\xi<\delta}\in M$. By elementarity and
$\crit(j)=\kappa$, $j(\seqof{\kappa_\xi}{\xi<\delta})=\seqof{\kappa_\xi}{\xi<\delta}$. 
Hence
\begin{xitemize}
\item[] 
  $N\modelof{j(\kappa)=\lim( j(\seqof{\kappa_\xi}{\xi<\delta}))=\lim(\seqof{\kappa_\xi}{\xi<\delta})=\kappa}$.
\end{xitemize}
This is contradiction to $\kappa=\crit(j)$ \footnote{Actually, we do not need the $\nu$-cc or any 
  other condition on $\calP$ to prove \assertof{b}. }. 

\assertof{c}\quad Finally, we prove that $\kappa$ is weakly Mahlo. Suppose that 
$C\subseteq\kappa$ is a club. Let $A\subseteq\kappa$ be \st\ it codes $C$ as well as 
witnesses of singularity of all singular cardinals and successorship of the successor 
cardinals $<\kappa$.

Let $M$ be a transitive model of $\ZFC^-$ \st\ $\kappa$, $A\in M$ and there is a $\nu$-cc \po\
$\poP$ with $(\uniV,\poP)$-generic $\genG$ \st\ there are $j$, $N\in\uniV[\genG]$ \st\
\xitemof{x-ccc-gen-0} holds.

Since $C\in M$ by $A\in M$ and $M\modelof{C\mbox{ is a club subset of }j(\kappa)}$, we have
$N\modelof{j(C)\mbox{ is a club subset of }\kappa}$ by elementarity. Since 
$j(C)\cap\kappa=C$ by $\crit(j)=\kappa$, it follows that $\kappa\in j(C)$. $\kappa$ is 
regular by \assertof{b}. Since $\poP$ preserves cardinality 
and cofinality $\GE\nu$ by its $\nu$-cc, $\uniV[\genG]\modelof{\kappa\mbox{ is regular}}$. It 
follows that $N\modelof{\kappa\mbox{ is regular}}$. Thus
$N\modelof{j(C)\xmbox{ contains a regular cardinal}}$ and 
$M\modelof{C\mbox{ contains a regular cardinal}}$ by elementarity.
By the choice of $A$ the weakly inaccessible cardinal in $C\cap M$ is really weakly 
inaccessible.

Since $C$ was arbitrary, this shows that $\kappa$ is a weakly Mahlo cardinal.
\smallskip

\assertof{2}: Suppose that $T$ is a $\kappa$-tree. We want to show that $T$ has 
a $\kappa$-branch. 

Since we have $\cardof{T}=\kappa$, we may assume \wolog\ that the underlying set of $T$ 
is $\kappa$.

Let $A\subseteq\kappa$ code the tree ordering $\leq_T$ as well as the witnesses asserting that $T$ is a
$\kappa$-tree. Let $M$ be a transitive model of $\ZFC^-$ 
\st\ $\kappa$, $A\in M$ and there is a $\nu$-cc \po\ 
$\poP$ with $(\uniV,\poP)$-generic $\genG$ \st\ there are $j$, $N\in\uniV[\genG]$ \st\
\xitemof{x-ccc-gen-0} holds.

We have $T\in M$ and $M\modelof{T\mbox{ is a $\kappa$-tree}}$ by $A\in M$. It follows 
that $N\modelof{j(T)\mbox{ is a $j(\kappa)$-tree}}$ by elementarity,  
and $j(T)_{\LT\kappa}=T$. Since $j(\kappa)>\kappa$, there is $t^*\in j(T)$  \st\
$N\modelof{t^*\in{j(T)}_\kappa}$. 

Let $\utilde{\leq}$ and $\utilde{t}$ be $\poP$-names 
of $j(\leq_T)$ and $t^*$. 

Back in $\uniV$,  let
\begin{xitemize}
\item[] $T_0:=\setof{t\in T}{\forces{\poP}{\check{t}\ \utilde{\leq}\ \utilde{t}}}$. 
\end{xitemize}
$T_0$ is a tree of height $\kappa$ and, by the $\nu$-cc of $\poP$, it is of width $\leq\nu$ 
and $\nu^+<\kappa$ (in 
$\uniV$). By a theorem of Kurepa (Proposition 7.90 in \cite{higher-inf}), it follows that 
there is a $\kappa$-branch $b$ in $T_0$. Clearly $b_0$ is also a $\kappa$-branch of $T$. 
\qedofLemma
\qedskip


\section{Generically measurable cardinals}
\Label{meas}

Let us call a cardinal $\kappa$ {\It greatly weakly Mahlo} if $\kappa$ is weakly inaccessible and 
there exists a non-trivial $\LT\kappa$-complete normal filter $\calF$ over $\kappa$ \st\
$\setof{\mu<\kappa}{\mu\xmbox{ is a regular cardinal}}\in\calF$, and $\calF$ is closed 
\wrt\ the Mahlo operation \footnote{Closedness here means that for any $S\in\calF$, we have $\mahlo(S)\in\calF$.}:
\begin{xitemize}
\xitem[x-gen-a-a-a-0]
  $S\ \mapsto\ {\darkred\mahlo(S)}:=\setof{\alpha\in S}{{}
  \begin{array}[t]{@{}l}
    \alpha\mbox{ has uncountable cofinality and}\\
    S\cap\alpha\mbox{ is stationary in }\alpha}.\end{array}$
\end{xitemize}
This definition of the Mahlo operation is slightly different from the one given in 
\cite{baumgartner-taylor-wagon}. 

For $\alpha\in\On$, we define the notion of $\alpha$-weakly Mahloness for all cardinals
$\kappa$ by induction on $\alpha$.
\begin{xitemize}
\xitem[x-gen-a-a-0] $\kappa$ is 
{\It $0$-weakly Mahlo} if $\kappa$ is weakly Mahlo;
\xitem[x-gen-a-0] $\kappa$ is 
{\It $1$-weakly Mahlo} if $\kappa$ is weakly Mahlo\footnotemark\ and
  $\setof{\mu<\kappa}{\mu\mbox{ is weakly Mahlo}}$\\ is stationary;
\xitem[x-gen-a-1] 
 for $1<\alpha\leq\kappa$,
  $\kappa$ is {\It$\alpha$-weakly Mahlo} if $\setof{\mu<\kappa}{\mu\mbox{ is $\beta$-weakly Mahlo}}$ is 
  stationary in $\kappa$ for all $\beta<\alpha$. 
\xitem[x-gen-a-2] 
  $\kappa$ is {\It hyper-weakly Mahlo} if $\triangle_{\alpha<\kappa}\setof{\mu<\kappa}{\mu
  \mbox{ is }\alpha\mbox{-weakly Mahlo}}$ is stationary\footnote{The definition of ``hyper 
  Mahloness'' (i.e. the strongly Mahlo version of the hyper-weakly Mahloness defined here) 
  has several deviations: in some cases $\kappa$-Mahloness (which 
  is apparently slightly weaker than the hyper Mahloness parallel to the hyper-weakly 
  Mahloness as defined here) is called hyper Mahlo.}. 
\end{xitemize}
\footnotetext{By \Lemmaof{P-gen-0-0}, the weak Mahloness of $\kappa$ follows from the second condition.}
\begin{Lemma}
  \Label{P-gen-0-0} For an ordinal $\kappa$, if $S\subseteq\kappa$ is a stationary set 
  consisting of regular cardinals, then $\kappa$ is also regular and hence $\kappa$ is 
  weakly Mahlo.
\end{Lemma}
\prf Suppose that $S$ is as above but $\kappa$ is not regular. 

We have $\cf{\kappa}>\omega$, since if $\cf(\kappa)=\omega$, then any 
increasing $\omega$-sequence of successor ordinals cofinal in $\kappa$ is a club 
in $\kappa$ disjoint from $S$. 

Say, $\cf(\kappa)=\mu<\kappa$. Let $\seqof{\xi_\alpha}{\alpha<\mu}$ be a continuously increasing 
sequence of ordinals cofinal in $\kappa$ \st\ $\xi_0>\mu$. By the assumption on $S$, there 
is $\lambda\in S\cap\setof{\xi_\alpha}{\alpha<\mu}$. Say, $\lambda=\xi_{\alpha^*}$. Then
$\cf(\lambda)\leq\alpha^*<\mu<\lambda$. This is a contradiction since $\lambda$ as an 
element of $S$ must be regular. 
\qedofLemma

\begin{Lemma}
  \Label{P-gen-0-0-0} Suppose $\alpha\leq\beta\leq\kappa$. If $\kappa$ is $\beta$-weakly 
  Mahlo, then $\kappa$ is $\alpha$-weakly Mahlo. 
\end{Lemma}
\prf
By induction on $\beta$. 
\qedofLemma
\qedskip

For $S\subseteq\kappa$ and $\alpha<\kappa$, let {\darkred$\mahlo^\alpha(S)$} be defined inductively by
\begin{xitemize}
\xitem[x-gen-0] 
  $\mahlo^0(S):=S$;
\xitem[x-gen-1] $\mahlo^{\alpha+1}(S):=\mahlo(\mahlo^\alpha(S))$;
\xitem[x-gen-2] 
  $\mahlo^\gamma(S):=\bigcap_{\alpha<\gamma}\mahlo^\alpha(S)$ for a limit $\gamma<\kappa$. 
\end{xitemize}
Finally, let
\begin{xitemize}
\xitem[x-gen-3] $\mahlo^\kappa(S):=\triangle_{\alpha<\kappa}\mahlo^\alpha(S)$. 
\end{xitemize}

Note that stationary sets are not necessarily closed \wrt\ intersection of decreasing 
sequence of short length: Let $\kappa$ be an uncountable cardinal with
$\kappa\geq\omega_\omega$. For $n\in\omega$, let
$S_n:=\setof{\alpha<\kappa}{\omega_n\leq\cf(\alpha)<\omega_\omega}$. Then each $S_n$, 
$n\in\omega$ is stationary. 
But $\bigcap_{n\in\omega}S_n=\emptyset$. 

\begin{Lemma}
  \Label{P-gen-1}\wassertof{1}\quad For a regular $\kappa$, a filter $\calF$ over $\kappa$ 
  is uniform (i.e. every end-segment of $\kappa$ is in $\calF$) and normal, if and only if 
  $\calF$ is non-principal, $\LT\kappa$-complete and normal.\smallskip

  \wassert{2} 
  If $\calF$ is a uniform normal filter over a regular $\kappa$, then $C\in\calF$ for all 
  club $C\subseteq\kappa$. It follows that all $S\in\calF$ are  
  stationary in $\kappa$. \smallskip 

  \wassert{3} If $\kappa$ is greatly weakly Mahlo and $\calF$ is as in the definition of 
  the greatly weak Mahloness of $\kappa$, then for all $\alpha<\kappa$
  $\setof{\xi<\kappa}{\xi\mbox{ is }\alpha\mbox{-Mahlo}}\in\calF$. 
\end{Lemma}
\prf \assertof{1}: ``$\Leftarrow$'' is trivial. For ``$\Rightarrow$'', suppose that 
$\delta<\kappa$ and $S_\alpha\in\calF$ for all $\alpha<\delta$.

For $\alpha<\kappa$, let 
\begin{xitemize}
\item[] $S^*_\alpha=\left\{\,
  \begin{array}{@{}ll}
    S_\alpha\setminus\delta, &\mbox{if }\alpha<\delta;\\
    \kappa &\mbox{otherwise}.
  \end{array}
  \right.$
\end{xitemize}
We have $S^*_\alpha\in\calF$ for $\alpha<\delta$ as $\kappa\setminus\delta\in\calF$ since 
$\calF$ is uniform. 

Then 
$\calF\ni\triangle_{\alpha<\kappa}S^*_\alpha=\bigcap_{\alpha<\kappa}S_\alpha\setminus\delta\subseteq 
\bigcap_{\alpha<\kappa}S_\alpha$.
Thus $\cap_{\alpha<\kappa}S_\alpha\in\calG$. \smallskip

\assertof{2}: We show first that
$Lim(\kappa)\ (=\setof{\alpha<\kappa}{\alpha\mbox{ is a limit ordinal}})$ is an element of
$\calF$. This follows from $Lim(\kappa)=\triangle_{\alpha<\kappa}\kappa\setminus(\alpha+1)\in\calF$. 

For a club $C\subseteq\kappa$, let $\seqof{c_\alpha}{\alpha<\kappa}$ be an increasing 
enumeration of $C$. Then we have
$C\supseteq Lim(\kappa)\cap\triangle_{\alpha<\kappa}\kappa\setminus c_\alpha\in\calF$. 

For an $S\in\calF$, $S\cap C\in\calF$ and hence $S\cap C\not=\emptyset$ for all club
$C\subseteq\kappa$. Thus $S$ is stationary in $\kappa$.\smallskip

\assertof{3}: By induction on $\alpha$. 
\qedofLemmaA

The following Proposition is a variant of Proposition 16.8 in \cite{higher-inf}. 

\begin{Prop}
  \Label{P-gen-0-1} Suppose that $\kappa$ is greatly weakly Mahlo, and let $\calF$ be a non-trivial
  $\LT\kappa$-complete normal filter over $\kappa$ \st\
  \begin{xitemize}
  \xitem[x-gen-4-0] $\Reg(\kappa):=\setof{\mu<\kappa}{\mu\mbox{ is regular\,}}\in\calF$, and
  \xitem[x-gen-4-1] $\calF$ is closed \wrt\ the Mahlo operation (as defined in 
    \xitemof{x-gen-a-a-a-0}). 
  \end{xitemize}
  Then, for any $1\leq\alpha<\kappa$,\smallskip

  \wassert{1} $\mahlo^\alpha(\Reg(\kappa))\in\calF$,\smallskip

  \wassert{2} $\mahlo^\alpha(\Reg(\kappa))=\setof{\mu<\kappa}{\mu\mbox{ is }
    \beta\mbox{-weakly Mahlo for all }\beta<\alpha}\\
  =\setof{\mu<\kappa}{\mu\mbox{ is }
    \alpha_0\mbox{-weakly Mahlo}}$ for all $1\leq\alpha<\omega$ where $\alpha_0$ is \st\ $\alpha=\alpha_0+1$;\\
  $\mahlo^\alpha(\Reg(\kappa))=\setof{\mu<\kappa}{\mu\mbox{ is }
    \beta\mbox{-weakly Mahlo for all }\beta\leq\alpha}\\
  =\setof{\mu<\kappa}{\mu\mbox{ is }
    \alpha\mbox{-weakly Mahlo}}$ for 
  all $\omega\leq\alpha<\kappa$,\smallskip

  \wassert{3} $\kappa$ is hyper-weakly Mahlo. 
\end{Prop}
\prf We first prove \assertof{1} and \assertof{2} simultaneously by induction on $1\leq\alpha<\kappa$.
Note that the last equality in both of the cases in \assertof{2} follows from 
\Lemmaof{P-gen-0-0-0}. 

For
$\alpha=1$, we have
\begin{xitemize}
\item[] $\calF\smash{\ubecause{\ni}{}{by \xitemof{x-gen-4-0} and \xitemof{x-gen-4-1}}}\mahlo(\Reg(\kappa))\,
  \begin{array}[t]{@{}l}
    =\mahlo^1(\Reg(\kappa))\\[\jot]
    =\setof{\mu\in\Reg(\kappa)}{\mu\cap\Reg(\kappa)\mbox{ is stationary in }\mu}\\[\jot]
    =\setof{\mu\in\Reg(\kappa)}{\mu\mbox{ is weakly Mahlo}}.\\[\jot]
    =\setof{\mu\in\Reg(\kappa)}{\mu\mbox{ is $0$-weakly Mahlo}}.
  \end{array}
  $
\end{xitemize}

Suppose that $\gamma<\kappa$ is a limit ordinal, and \assertof{1}, \assertof{2} hold for 
all $\alpha<\gamma$. Then

\begin{xitemize}
\item[] $\mahlo^\gamma(\Reg(\kappa))\ubecause{=}{}{by \xitemof{x-gen-2}}
\bigcap_{\alpha<\gamma}\mahlo^\alpha(\Reg(\kappa))\in\calF$
\end{xitemize}
by the induction hypothesis about \assertof{1} and $\LT\kappa$-completeness of $\calF$. 

Suppose that $\mu\in\mahlo^\gamma(Reg(\kappa))$. Then, by the induction hypothesis about 
\assertof{2}, $\mu$ is $(\beta+1)$-weakly Mahlo for all $\beta<\gamma$. By 
\xitemof{x-gen-a-1}, it follows that
$\setof{\xi<\mu}{\xi\mbox{ is }\beta\mbox{-weakly Mahlo}}$ is stationary in $\mu$. Thus, 
again by \xitemof{x-gen-a-1}, $\mu$ is $\gamma$-weakly Mahlo.

Conversely, if $\mu<\kappa$ is $\gamma$-weakly Mahlo, then, by \Lemmaof{P-gen-0-0-0}, $\mu$ 
is $\alpha$-weakly Mahlo for all $\alpha<\gamma$. Thus, by the induction hypothesis about 
\assertof{2},\\
$\mu\in\bigcap_{\alpha<\gamma}\mahlo^\alpha(\Reg(\kappa))=\mahlo^\gamma(\Reg(\kappa))$.  
This shows that \assertof{2} holds for $\gamma$. \smallskip

Suppose now that \assertof{1} and \assertof{2} hold for  $1\leq\alpha<\kappa$.

If $\alpha<\omega$, this means in particular that for $\alpha_0$ \st\ $\alpha=\alpha_0+1$,
\begin{xitemize}
\item[] 
  $\mahlo^\alpha(\Reg(\kappa))=\setof{\mu<\kappa}{\mu\mbox{ is }\alpha_0\mbox{-weakly Mahlo}}\in\calF$. 
\end{xitemize}

By the definition \xitemof{x-gen-1} of the iteration of Mahlo operation and 
\xitemof{x-gen-4-1}, we have 
\begin{xitemize}
\item[] 
  $\mahlo^{\alpha+1}(\Reg(\kappa))
  =\setof{\mu<\kappa}{\mu\mbox{ is }\ubecause{(\alpha_0+1)}{}{$=\alpha$ }\mbox{-weakly Mahlo}}\in\calF$.
\end{xitemize}

If $\omega\leq\alpha<\omega$, our assumption is 
\begin{xitemize}
\item[] 
  $\mahlo^\alpha(\Reg(\kappa))=\setof{\mu<\kappa}{\mu\mbox{ is }\alpha\mbox{-weakly Mahlo}}\in\calF$. 
\end{xitemize}
Thus, similarly to above, we obtain
\begin{xitemize}
\item[] 
  $\mahlo^{\alpha+1}(\Reg(\kappa))
  =\setof{\mu<\kappa}{\mu\mbox{ is }(\alpha+1)\mbox{-weakly Mahlo}}\in\calF$.
\end{xitemize}

\assertof{1} and \assertof{2} imply \assertof{3}:  
\begin{xitemize}
\item[] 
  $\triangle_{\alpha<\kappa}\setof{\mu<\kappa}{\mu
  \mbox{ is }\alpha\mbox{-weakly Mahlo}}\ubecause{=}{}{by \assertof{2}}\triangle_{\alpha<\kappa}
  \mahlo^\alpha(\Reg(\kappa))
  \ubecause{\in}{}{by \assertof{1} and normality of $\calF$}\calF$. 
\end{xitemize}
In particular $\triangle_{\alpha<\kappa}\setof{\mu<\kappa}{\mu
  \mbox{ is }\alpha\mbox{-weakly Mahlo}}$ is stationary, and this proves that $\kappa$ is 
hyper-weakly Mahlo.
\qedofProp
\qedskip

The following theorem actually holds already for a $\nu$-cc-g. weakly compact $\kappa$ for 
some $\nu<\kappa$. This will be addressed in the forthcoming \cite{fuchino-sakai-2}. 

\begin{Thm}
  \Label{P-gen-2} If $\kappa$ is a $\nu$-cc-g.\ measurable cardinal for a $\nu<\kappa$, 
  then $\kappa$ is greatly weakly Mahlo. 
\end{Thm}
\prf Let $\poP$ be a ccc \po\ with $(\uniV,\poP)$-generic $\genG$ \st\ there are 
classes $j$, $M\subseteq\uniV[\genG]$ with $\Elembed{j}{\uniV}{M}{\kappa}$.

Note that, since generically large cardinals are definable (see \cite{fuchino-sakai}), we 
may apply forcing theorems in the arguments which involve $j$ and $M$. In particular, we 
may assume that
\begin{xitemize}
\xitem[x-gen-5] 
  $\forces{\poP}{\Elembed{j}{\uniV}{M}{\kappa}}$. 
\end{xitemize}

In $\uniV[\genG]$, let $\tilde{\calF}:=\setof{S\subseteq\kappa}{S\in\uniV, j(S)\ni\kappa}$ 
and let $\utilde{\calF}$ be a $\poP$-name of $\tilde{\calF}$.

In $\uniV$, 
let $\calF:=\setof{S\subseteq\kappa}{\forces{\poP}{\check{S}\varin\utilde{\calF}}}
=\setof{S\subseteq\kappa}{\forces{\poP}{j(\check{S})\ni\kappa}}$. Then
\begin{Claim}
  \Label{cl-gen-0}
  \wassertof{1} $\calF$ is a non-trivial $\LT\kappa$-complete normal filter.\smallskip

  \wassert{2} $\Reg(\kappa)\in \calF$.\smallskip

  \wassert{3} $\calF$ is closed \wrt\ Mahlo operation. 
\end{Claim}
\prfofClaim 
\assertof{1}: It is clear that $\calF$ is a non-trivial filter. 

Suppose that $\vec{S}:=\seqof{S_\alpha}{\alpha<\mu}\in\uniV$ for some $\mu<\kappa$ 
is a sequence of length $\mu$ of elements of $\calF$. Then
$\forces{\poP}{\vec{S}\xmbox{ is a sequence of elements of $\utilde{\calF}$ of length $\mu$}}$.
Since 
$\forces{\poP}{j(\vec{S})=\seqof{j(S_\alpha)}{\alpha<\mu}}$ by
\xitemof{x-gen-5},  
we have
\begin{xitemize}
\item[] $\forces{\poP}{j(\bigcap\vec{S})=\bigcap j(\vec{S})
  =\bigcap\setof{j(S_\alpha)}{\alpha<\mu}\ni\kappa}$.
\end{xitemize}
Thus $\forces{\poP}{\bigcap\vec{S}\in\utilde{\calF}}$, and hence 
$\bigcap\vec{S}\in\calF$.

If $\vec{S}:=\seqof{S_\alpha}{\alpha<\kappa}$ is a sequence in $\uniV$ of elements of $\calF$. 
then 
\begin{xitemize}
\item[] 
  $\forces{\poP}{\kappa\in\bigcap\setof{j(S_\alpha)}{\alpha<\kappa}
    =\bigcap\big((j(\vec{S}))\restr\kappa}\big)$. 
\end{xitemize}
Since 
  $\forces{\poP}{\kappa\in\bigcap\big((j(\vec{S}))\restr\kappa\big)\ \Leftrightarrow\ 
  \kappa\in\triangle j(\vec{S})^\bullet=j(\triangle\vec{S})}$,
it follows that $\forces{\poP}{\triangle\vec{S}\in\utilde{\calF}}$ and thus 
$\triangle\vec{S}\in\calF$. \smallskip

\assertof{2}: Let $R:=\Reg(\kappa)\setminus\nu$. Then we have
$\forces{\poP}{j(R)=\Reg(j(\kappa))^M\setminus\nu}$ by \xitemof{x-gen-5}. 
By the $\nu$-cc of $\poP$, it follows that
$\forces{\poP}{\kappa\mbox{ is regular and  }\kappa>\nu}$. Thus 
$\forces{\poP}{\kappa\in j(R)}$, and hence $R\in\calF$. \smallskip

\assertof{3}: If $S\in\calF$, 
then $S$ is stationary by \assertof{1} and \Lemmaof{P-gen-1},\,\assertof{2}. 

Since $\poP$ 
is $\nu$-cc, $\forces{\poP}{V^\poP\models S\xmbox{ is stationary in }\kappa}$. 

Since $\forces{\poP}{S=j(S)\cap\kappa}$, 
it follows that 
\begin{xitemize}
\item[] 
  $\forces{\poP}{\kappa\in\mahlo^M(j(S))=j(\mahlo^\uniV(S))}$.  
\end{xitemize}
Thus, $\forces{\poP}{\mahlo^\uniV(S)\in\utilde{\calF}}$ and hence $\mahlo(S)\in\calF$. 
\qedofClaim\\
\qedofThm
\qedskip

\begin{Prop}
  \Label{P-gen-3} For a regular cardinals $\kappa$, $\nu$ with $\nu<\kappa$, 
  \tfae:\smallskip

  \wassert{a} $\kappa$ is $\nu$-cc-g.\ measurable.

  \wassert{b} There is a non-trivial, non-principal and $\nu$-saturated $\LT\kappa$-complete ideal over $\kappa$.

  \wassert{c} there are a $\nu$-cc \po\/ $\poP$, a $(\uniV,\poP)$-generic filter $\genG$, 
  and $j$, $M\subseteq\uniV[\genG]$ \st\/
  $\uniV[\genG]\modelof{\Elembed{j}{\uniV}{M}{\kappa}}$ and
  $(\fnsp{\kappa}{M})^{\uniV[\genG]}\subseteq M$.
\end{Prop}
\prf ``\assertof{c} $\Rightarrow$ \assertof{a}'': is clear. So we shall prove ``\assertof{a}
$\Rightarrow$ \assertof{b}'' and ``\assertof{b} $\Rightarrow$ \assertof{c}''.
\smallskip

``\assertof{a} $\Rightarrow$ \assertof{b}'': Let $\poP$ be a $\nu$-cc \po\ \st, for $(\uniV,\poP)$-generic 
$\genG$ and $j$, $M\subseteq\uniV[\genG]$, we have $\Elembed{j}{V}{M}{\kappa}$.

In $\uniV$, let $\calI:=\setof{A\subseteq\kappa}{\forces{\poP}{\kappa\not\in j(\check{A})}}$. 
Note that $\calI$ is the dual ideal of the filter of $\calF$ in the proof of 
\Propof{P-gen-0-1}.

\begin{Claim}
  $\calI$ is $\LT\kappa$-complete and $\nu$-saturated ideal (in $\uniV$).
\end{Claim}
\prf $\LT\kappa$-completeness follows from \Claimof{cl-gen-0},\,\assertof{1}.

In the following, we argue in $\uniV$. 
To prove that $\calI$ is $\nu$-saturated, assume, toward a contradiction, that
$\seqof{A_\xi}{\xi<\nu}$  is a pairwise incompatible sequence of elements in
$\poP(\kappa)\setminus\calI$. By the $\LT\kappa$-completeness of $\calI$, we may choose the sequence \st\
$A_\xi$, $\xi<\nu$ are pairwise disjoint. For each $\xi<\nu$, since $A_\xi\not\in\calI$, there 
is $\condp_\xi\in\poP$, \st\ $\condp_\xi\forces{\poP}{\kappa\in j(\check{A_\xi})}$. 
By $\nu$-cc of $\poP$, there are $\xi<\eta<\nu$ \st\ $\condp_\xi$ and $\condp_\eta$ are 
compatible, say $\condr\leq_\poP \condp_\xi$, $\condp_\eta$. But then
$\condr\forces{\poP}{\kappa\in j(\check{A}_\xi)\cap j(\check{A_\eta})}$ and hence
$\condr\forces{\poP}{\check{A}_\xi\cap \check{A_\eta}\not=\emptyset}$. It follows that
$A_\xi\cap A_\eta\not=\emptyset$. This is a contradiction to the choice 
of $A_\xi$, $\xi<\nu$. 
\qedofClaim
\qedskip

``\assertof{b} $\Rightarrow$ \assertof{c}'': Let $\calI$ be 
a $\nu$-saturated $\kappa$-complete ideal over $\kappa$. $\poP_\calI:=\poP(\kappa)\setminus\calI$ 
satisfies then the $\nu$-cc. 

$\calI$ is precipitous by Lemma 22.22 in \cite{millennium-book}. 
Let $\genG$ be a $(\uniV,\poP_\calI)$-generic filter, and let $\Elembed{j}{\uniV}{M}{\kappa}$ be 
the canonical elementary embedding of $\uniV$ into the Mostowski collapse of the generic 
ultrapower by $\genG$. By Lemma 22.31 in \cite{millennium-book}, we have
$(\fnsp{\kappa}{M})^{\uniV[\genG]}\subseteq M$. \qedofProp

\begin{Thm}
  \Label{P-gen-4} Suppose that $\kappa$ is a $\nu$-cc-g.\ measurable cardinal for some 
  regular $\nu<\kappa$. 
  Then $\kappa$ is the stationary limit 
  of $\nu$-cc-g.\ weakly compact cardinals.
\end{Thm}
\prf Suppose that $\kappa$ is $\nu$-cc-g.\ measurable and let $C\subseteq\kappa$ be an 
arbitrary club subset of $\kappa$. We have to show that $C$ contains a $\nu$-cc-g.\ weakly 
compact cardinal. Let $\poP$ be a $\nu$-cc \po\ with a $(\uniV,\poP)$-generic
$\genG$ and $j$, $M\subseteq\uniV[\genG]$ \st\ $\Elembed{j}{\uniV}{M}{\kappa}$ and 
\begin{xitemize}
\xitem[x-gen-6] $(\fnsp{\kappa}{M})^{\uniV[\genG]}\subseteq M$
\end{xitemize}
(see \Propof{P-gen-3}). 

Since $M\modelof{j(C)\mbox{ is a club subset of }j(\kappa)}$ and $\kappa\in j(C)$ by the closedness, the 
following claim completes the proof.
\begin{Claim}
  \Label{cl-gen-1} $M\modelof{\kappa\mbox{ is }\nu\mbox{-cc-g.\ weakly compact\/}}$. 
\end{Claim}
\prfofClaim In $M$, suppose $A\subseteq\kappa$. We have to show in $M$ that there is a 
transitive model $M_0$ of $\ZFC^-$ with $\kappa$, $A\in M_0$ and 
$\Elembed{j_0}{M_0}{N_0}{\kappa}$ for some $j_0$, $N_0$ in some $\nu$-cc generic extension.

By \Propof{P-gen-3}, there is a $\nu$-saturated $\LT\kappa$-complete ideal $\calI$ on
$\kappa$ in $\uniV$.

In $\uniV[\genG]$, let $\calJ$ be the ideal over $\kappa$ generated by $\calI$. By $\nu$-cc 
of $\poP$, it is easy 
to see that $\calJ$ is $\LT\kappa$-complete (in $\uniV[\genG]$).  $\calJ$ is $\nu$-saturated (in
$\uniV[\genG]$) by Prikry's Theorem (see e.g.\ Theorem 17.1 in \cite{higher-inf}).  
Let $\poP_\calJ:=(\psof{\kappa}\setminus\calJ)^{\uniV[\genG]}$. Then
$\uniV[\genG]\modelof{\poP_\calJ\mbox{ has the }\nu\mbox{-cc}}$. 

Working further in $\uniV[\genG]$, let $\theta$ be sufficiently large and let $M_1$ be \st\
\begin{xitemize}
\xitem[x-gen-7] $M_1\prec\calH(\theta)$,
\xitem[x-gen-8] $\kappa+1\cup\ssetof{A,\calJ}\subseteq M_1$, and 
\xitem[x-gen-10] $\cardof{M_1}=\kappa$. 
\end{xitemize}

Let $\mapping{m}{M_1}{M_0}$ be the Mostowski collapse. Note that 
\begin{xitemize}
\xitem[x-gen-11] $m\restr\kappa+1=\id_{\kappa+1}$
\end{xitemize}
by \xitemof{x-gen-8}. $M_0\in M$ by \xitemof{x-gen-6}. 


By $A\subseteq\kappa$ and \xitemof{x-gen-11}, we have $A=m(A)\in M_0$. 

Let $\calJ_0:=m(\calJ)$. $\calJ_0\in M_0$ by this definition and $\calJ_0=\calJ\cap M_0$ by 
\xitemof{x-gen-11}. 
By the elementarity \xitemof{x-gen-7} (and since $\theta$ is taken sufficiently large), we 
have
$M_1\modelof{\calJ\mbox{ is a }\nu\mbox{-saturated, }\LT\kappa\mbox{-complete ideal over }\kappa}$. 
It follows that 
\begin{xitemize}
\item[] 
  $M_0\modelof{\calJ_0\mbox{ is a }\nu\mbox{-saturated, }\LT\kappa\mbox{-complete ideal over	}\kappa}$.
\end{xitemize}
In particular, $\calJ_0$ is precipitous in connection with $M_0$ by Lemma 22.22 in 
\cite{millennium-book}. 

Let $\poQ:=(\psof{\kappa}\setminus\calJ_0)^{M_0}$ (note that $\poQ\in M$ since $M_0\in M$). 
$m^{-1}\restr\poQ=\id_\poQ$ is then an order-preserving and incompatibility preserving embedding of $\poQ$ 
into $\poP_\calJ$. Since $\poP_\calJ$ is $\nu$-cc in $\uniV[\genG]$, $\poQ$ is 
also $\nu$-cc in $\uniV[\genG]$. It follows that $\poQ$ is also $\nu$-cc in $M$ (note that
$(\nu^+)^M=(\nu^+)^{\uniV[\genG]}$ by \xitemof{x-gen-6}). 

Let $\genH$ be a $(M,\poQ)$-generic filter and let $\Elembed{j_0}{M_0}{N_0}{\kappa}$ where 
$N_0$ is the Mostowski collapse of the generic ultrapower of $M_0$ by $\genH$ 
(in $M[\genH]$).

Clearly,  $M_0$ together with these $j_0$ and $N_0$ is as desired. 
\qedofClaim\\
\qedofThm
\section{Generically supercompact cardinals}
\Label{gen-sc}

\begin{Thm}
  \Label{P-gen-5} Suppose that $\kappa$ is $\nu$-cc-g.\ $2^\kappa$-supercompact for some 
  uncountable cardinal $\nu<\kappa$. Then $\kappa$ is the stationary limit of $\nu$-cc-g.\ 
  measurable cardinals.
\end{Thm}
\prf Let $\poP$ be a $\nu$-cc \po\ with a $(\uniV,\poP)$-generic filter $\genG$ and $j$,
$M\subseteq\uniV[\genG]$ \st, in $\uniV[\genG]$, $\Elembed{j}{\uniV}{M}{\kappa}$,
$j(\kappa)>(2^\kappa)^\uniV$ and
\begin{xitemize}
\xitem[x-gen-13] 
  $j\imageof(2^\kappa)^\uniV\in M$. 
\end{xitemize}

As in the proof of \Thmof{P-gen-4}, it is enough to show that
$M\modelof{\kappa\mbox{ is }\nu\mbox{-cc-g.}\xmbox{ measurable}}$. Thus, by \Propof{P-gen-3}, we are 
done by showing that
$M\modelof{\mbox{there is a }\nu\mbox{-saturated, }\LT\kappa\xmbox{-complete ideal over }\kappa}$. 

Since $\kappa$ is $\nu$-cc g.\ measurable, there is a $\nu$-saturated, $\LT\kappa$-complete 
ideal $\calI$ over $\kappa$ in $\uniV$ by \Propof{P-gen-3}. In $\uniV[\genG]$, let $\calJ$ 
be the ideal over $\kappa$ generated by $\calI$. By Prikry's theorem, we have
\begin{xitemize}
\xitem[x-gen-14] 
  $\uniV[\genG]\modelof{\calJ\mbox{ is }\nu\mbox{-saturated, }\LT\kappa\mbox{-complete ideal}}$.
\end{xitemize}

Note that
\begin{xitemize}
\item[] $
  \begin{array}{@{}r@{{}={}}l}
    \calJ &\setof{A\in\psof{\kappa}^{\uniV[\genG]}}{A\subseteq B\mbox{ for some }B\in\calI}\\[\jot]
    &\setof{A\in\psof{\kappa}^{\uniV[\genG]}}{A\subseteq j(B)\mbox{ for some }B\in\calI}\\[\jot]
    &\setof{A\in\psof{\kappa}^{\uniV[\genG]}}{A\subseteq C\mbox{ for some }C\in j\imageof\calI}.
  \end{array}
  $
\end{xitemize}
Since $j\imageof\calI\in M$ by \xitemof{x-gen-13}, it follows that 
\begin{xitemize}
\item[] 
  $\calJ\cap M=\setof{A\in\psof{\kappa}^M}{A\subseteq C\mbox{ for some }C\in j\imageof\calI}$
\end{xitemize}
is an element of $M$. 

Thus, we have 
$M\modelof{\calJ\cap M\mbox{ is }\nu\mbox{-saturated, }\LT\kappa\mbox{-complete ideal}}$
by \xitemof{x-gen-14}. \qedofThm

\section{Reflection properties down to $\LT$ a generically supercompact cardinal}
\Label{refl}

An interesting fact about the notion of generic large cardinals is that the continuum can be 
generic large, or in some cases the continuum can be strictly larger than many generic 
large cardinals (cf.\ \Thmof{P-gen-4}, \Thmof{P-gen-5}). This is in particular the case 
with ccc-generically supercompact cardinals  
(e.g.\ obtained by starting with a supercompact $\kappa$ and then by forcing with $\Fn(\kappa,2)$). The largeness 
properties of generic supercompact cardinals for forcing with chain condition discussed in 
the previous sections can thus also be situations with the continuum. 

Since generic large cardinals are reflection points of diverse reflection statements (as discussed 
below), the continuum can be also the reflection point of the same reflection statements.

Let $\calS$ be a class of (not necessarily first-order) structures with a notion
$\sqsubseteq_\calS$ of substructure relation where $(\calS,\sqsubseteq_\calS)$ should 
satisfy certain reasonable properties like that 
\begin{xitemize}
\xitem[x-gen-15] 
  $\strA\in\calS$ and $\strA\cong\strB$ imply
  $\strB\in\calS$,\\[\jot]
  $\sqsubseteq_\calS$ is transitive, \\[\jot]
  $\strA\sqsubseteq_\calS\strB$ and $\pairof{\strB,\strA}\cong\pairof{\strB',\strA'}$ imply 
  $\strA'\sqsubseteq_\calS\strB'$,
\end{xitemize}
etc.\memo{ペーパーとして出版するときには，性質を全部書き下す．
}

We also assume that $\calS$ is absolute and $\sqsubseteq_\calS$ is upward absolute, meaning 
that if $M$, $N$ are transitive (class or set) models of $\ZFC^-$ and $M$ is an inner model 
of $N$, then, for any $\strA$, $\strB\in M$, $M\models\strA\in\calS$ $\Leftrightarrow$
$N\models\strA\in\calS$ and $M\models\strA\sqsubseteq_\calS\strB$ $\Rightarrow$
$N\models\strA\sqsubseteq_\calS\strB$. For a structure $\strA$, we denote by {\darkred$\ulsetof{\strA}$} 
the underlying set of $\strA$ and by {\darkred$\ccardof{\strA}$} the cardinality of the 
underlying set of the structure $\strA$. 

Given such a class $\calS=(\calS,\sqsubseteq_\calS)$ of structures, and a property $P$, the 
reflection number of $P$ (in connection with $\calS$) is defined as 
\begin{xitemize}
\item[]\mbox{}\hspace{-2em}
 ${\darkred\refl_{stat}(\calS,P)}:=\min\setof{\kappa\in\Reg}{{}
  \begin{array}[t]{@{}l}
    \mbox{for any }\strA\in\calS\mbox{ with }\strA\models P\mbox{ and }\ccardof{\strA}\geq\kappa,\\
    \mbox{the set }\setof{\ulsetof{\strB}}{\strB\sqsubseteq_\calS\strA,\strB\models 
      P,\,\ccardof{\strB}<\kappa}\\
  \mbox{is stationary subset of }[\,\ulsetof{\strA}\,]^{<\kappa}\quad}
  \end{array}
  $ 
\end{xitemize}
where we define $\min\emptyset:=\infty$.

The reflection spectrum of $P$ is 
\begin{xitemize}
\item[]\mbox{}\hspace{-2em} 
 ${\darkred\REFL_{stat}(\calS,P)}:=\setof{\kappa\in\Reg}{{}
  \begin{array}[t]{@{}l}
    \mbox{for any }\strA\in\calS\mbox{ with }\strA\models P\mbox{ and }\ccardof{\strA}\geq\kappa,\\
    \mbox{the set }\setof{\ulsetof{\strB}}{\strB\sqsubseteq_\calS\strA,\strB\models 
      P,\,\ccardof{\strB}<\kappa}\\
  \mbox{is stationary subset of }[\,\ulsetof{\strA}\,]^{<\kappa}\quad}
  \end{array}
  $ 
\end{xitemize}

\begin{Ex}
  \Label{Ex-gen-0} Let $\calS$ be the class of all first countable topological spaces $X=\pairof{X,\tau}$ 
  where $\tau$ is an open basis for the space. $\sqsubseteq_\calS$ is the subspace 
  relation.  

  For $P=$ non-metrizability, the consistency of $\refl_{stat}(\calS,P)=\aleph_2$ is known 
  as Hamburger's Problem which is open for almost a half century. 
\end{Ex}

A property $P$ is {\It downward absolute} if, for any transitive (class or set) models $M$, $N$ 
of $\ZFC^-$ \st\ $M$ is an inner model of $N$, and for any structure $\strA$, if
$N\modelof{\strA\models P}$ implies $M\modelof{\strA\models P}$.

For a class $\calP$ of \pos, {\It$\calP$ preserves $P$}, if $\strA\models P$ then for any
$\poP\in\calP$, $\forces{\poP}{\strA\models P}$ holds.

\begin{Thm}
  \Label{P-gen-6} Suppose that $\calS=(\calS,\sqsubseteq_\calS)$ is a class of structures, 
$\calP$ a class of \pos, and $\kappa$ a $\calP$-g.\ supercompact cardinal. If a 
  property $P$ satisfies: 
  \begin{xitemize}
  \xitem[x-gen-15-0] 
    $P$ is downward absolute and 
  \xitem[x-gen-15-1] $\calP$ preserves $P$, 
  \end{xitemize}
  then $\kappa\in\REFL_{stat}(\calS,P)$ and hence $\refl_{stat}(\calS,P)\leq\kappa$. 
\end{Thm}
\prf Suppose that $\strA\in\calS$ and $\ccardof{\strA}\geq\kappa$. By replacing $\strA$ 
with an isomorphic structure, we may assume that $\ulsetof{\strA}=\lambda\in\Card$ (see 
\xitemof{x-gen-15}).

Let $\poP\in\calP$ be \st\ for an $(\uniV,\calP)$-generic $\genG$ and $j$,
$M\subseteq\uniV[\genG]$, we have $\uniV[\genG]\modelof{\Elembed{j}{\uniV}{M}{\kappa}}$,
\begin{xitemize}
\xitem[x-gen-16] $j(\kappa)>\mu$ and $j\imageof{\mu}\in M$ where $\mu=\beth_n(\lambda)$ for 
  sufficiently large $n$.\footnotemark
\end{xitemize}
\footnotetext{We need $\beth_n(\lambda)$ here since the elements of $\calS$ may be ($n+1$)-th 
  order structures for $n\geq 1$. }
By the last condition \xitemof{x-gen-16}, we have $j(\strA)\restr j\imageof{\lambda}\in M$. 

$\uniV[\genG]\modelof{\strA\models P}$ by \xitemof{x-gen-15-1}.
By $\uniV[\genG]\modelof{\strA\cong j(\strA)\restr j\imageof{\lambda}}$, it follows that 
$\uniV[\genG]\modelof{j(\strA)\restr j\imageof{\lambda}\models P}$. 
Since $\ulsetof{j(\strA)\restr j\imageof{\lambda}}=j\imageof\lambda\in M$ by \xitemof{x-gen-16}, we 
have $j(\strA)\restr j\imageof{\lambda}\in M$ 
and hence $M\modelof{\strA\restr j\imageof{\lambda}\models P}$ 
by \xitemof{x-gen-15-0}.

Suppose $\calD$ is an arbitrary club subset of $[\lambda]^{\LT\kappa}$ in $\uniV$. 
Since $M\modelof{j\imageof\calD\xmbox{ is cofinal in }[j\imageof\lambda]^{\LT\kappa}}$,
we have $M\modelof{j\imageof\lambda=\bigcup(j\imageof\calD)\in j(\calD)}$. 

In $\uniV$, let $\calS_0=\setof{\strB}{\strB\sqsubseteq_\calS\strA, \ccardof{\strB}<\kappa}$. 
Since $j(\ulsetof{\strB})=j\imageof\ulsetof{\strB}$ for all $\strB\in\calS_0$, it follows 
that $\bigcup\setof{\ulsetof{\strC}}{\strC\in j\imageof\calS_0}=j\imageof\lambda$. 
Thus, $M\modelof{j(\strA)\restr j\imageof{\lambda}\sqsubseteq_\calS j(\strA)}$  
(for this, we have to assume that $\calS$ satisfies the the property that the union of 
upward directed system of $\sqsubseteq_\calS$-substructures is a $\sqsubseteq_S$-substructure 
and certain Downward L\"owenheim-Skolem theorem on $\sqsubseteq_\calS$-substructures of a 
given structure in $\calS$). \memo{これらの
  性質を書き出す．}
\wrt\ $\sqsubseteq_\calS$).

By elementarity, it follows that
$\setof{\ulsetof{\strB}}{\strB\sqsubseteq_\calS\strA,\,\ccardof{\strB}<\kappa}$ is 
stationary in $[\lambda]^{\LT\kappa}$. \qedofThm\qedskip

The following are some application of the theorem above.

\begin{Cor}
  Suppose that $\kappa$ is a ccc-g.\ supercompact cardinal and $\calS=(\calS,\leq)$ 
  is a variety and $\strA\in\calS$ with $\ccardof{\strA}\geq\kappa$ is not free. Then 
  there are  stationarily many non-free $\strB\leq\strA$ of cardinality $<\kappa$. In 
  particular, 
  with $P$ being ``being non-free'', we have $\REFL_{stat}(\calS, P)\ni\kappa$ and
  $\refl_{stat}(\calS, P)\leq\kappa$. 
\end{Cor}
\prf ccc \pos\ preserve non-freeness of algebras in any variety (see \cite{potential}).\\ \qedofCor
\qedskip

Note that the Corollary above applies e.g.\ to groups, abelian groups, Boolean algebras, etc.

We shall call \pos\ of the form $\Fn(\lambda,2)$ {\It generalized Cohen \pos}. For
$\calP=\setof{\poP}{\poP\mbox{ is forcing equivalent to some }\Fn(\lambda,2)}$,
$\calP$-g.\ supercompactness for this $\calP$ will be also called {\It Cohen-g.\ supercompactness}. 

In the following Corollaries,  stationarity may be replaced with clubness since any topological 
space containing a non-metrizable subspace is non-metrizable and any tree containing a 
non-special subtree is non-special. 

\begin{Cor} Suppose that $\kappa$ is a Cohen-g.\ supercompact cardinal. Then any 
  first countable non-metrizable topological space of cardinality $\geq\kappa$ have club 
  many subspaces of size $\LT\kappa$ which are non-metrizable. 
\end{Cor}
\prf Generalized Cohen \pos\ preserve non-metrizability (see Dow, Tall and Weiss 
\cite{dtw}). \qedofCor

\begin{Cor}
  Suppose that $\kappa$ is a Cohen-g.\ supercompact. Then any non-special tree $T$ has club 
  many non-special subtrees of size $<\kappa$. 
\end{Cor}
\prf Generalized Cohen \pos\ preserve non-specialty of trees (Todor\v{c}evi\'c, see 
\cite{RIMS15}). \qedofCor
\qedskip

The Cohen-g.\ supercompactness in \Corabove\ cannot be replaced by ccc-g.\ 
supercompactness: 

We can prove (in \ZFC) that there is a non-special tree of size $2^{\aleph_0}$ 
without branches of length $\omega_1$ (e.g.\
$T:=\setof{t}{\mapping{t}{\alpha}{\omega},\,t\mbox{ is 1-1 for some }\alpha<\omega_1}$ with 
the ordering $t/\leq_Tt'$ $:\Leftrightarrow$ $t\subseteq t'$ is such a tree).
\ifprivate{\privatecolor
Let $T$ be as above. It is clear that $T$ does not have any $\omega_1$ branch: such an 
branch $b$ in $T$ would give an injection $\mapping{\bigcup b}{\omega_1}{\omega}$.

To show that $T$ is non-special, assume, toward a contradiction,  
that $\mapping{f}{T}{\omega}$ is a specializing mapping. Note that $f$ restricted to any 
branch in $T$ is 1-1. 

Let $T_0:=\setof{t\in T}{\omega\setminus t\imageof{\dom(t)}\mbox{ is infinite}}$. 
For $s\in[\omega]^{<\aleph_0}$ let
\begin{xitemize}
\xitemA[x-gen-16-a] 
  $P_s(t):=\setof{n\in\omega}{{}
  \begin{array}[t]{@{}l}
    \mbox{there is }t'\in T_0\mbox{ \st\ }
    t\leq_T t',\ \\
    t'\imageof{\dom(t)}\cap s=\emptyset\mbox{, and }f(t')=n}.
  \end{array}$ 
\end{xitemize}

Let $\seqof{t_n}{n\in\omega}$ be a $\leq_T$-increasing sequence in $T_0$ and 
$\seqof{k_n}{n\in\omega}$ a 1-1 sequence with $k_n\in\omega$, $n\in\omega$ \st
\begin{xitemize}
\xitemA[] $t_0=\emptyset$;
\xitemA[]
  $k_n=\min(\omega\setminus(t_n\imageof{\dom(t_n)}\cup\ssetof{k_0\ctentenc k_{n-1}}))$ for
  $n\in\omega$;
\xitemA[x-gen-16-a-0]
  $t_{n+1}$ is \st\
  $\omega\setminus t\imageof{\dom(t_{n+1})}\supseteq\ssetof{k_0\ctenten  k_n}$;
\xitemA[x-gen-16-a-1] 
  If $P_{\ssetof{k_0\ctentenc k_n}}(t_n)\ni n$, then $f(t_{n+1})=n$.
\end{xitemize}

Let $t_\infty:=\bigcup_{n\in\omega}t_n$. Then $t_\infty\in T_0$ by \xitemof{x-gen-16-a-0} 
and thus, by \xitemof{x-gen-16-a-1}, no $n\in\omega$ can be assigned $t_\infty$ as
$f(t_\infty)$. This is a contradiction. 
}\fi 

By Baumgartner, Malitz, and Reinhardt \cite{BMR}, all trees of size $\LT\continuum$ without 
branches of length $\omega_1$ are special under Martin's Axiom.

If we start from a supercompact $\kappa$ and force Martin's Axiom together with
$\continuum=\kappa$ by the standard forcing construction, then in the resulting model
$\kappa$ $(=\continuum)$ is ccc-g.\ supercompact and MA holds. Hence, by the remark above, 
the statements in \Corabove\ does not hold in this model. If we denote the reflection 
number $\refl(\calC,F)$ for $\calC:=$ trees and $F:=$ non-special by $\refl_{\RC}$, the 
non-reflection mentioned above can be reformulated as: 
\MA\ implies $\refl_{\RC}>\continuum$. 

\ifprivate{\privatecolor
\medskip{}
If we force $\MA$ $+$ $\mu=\continuum$ for $\mu<\kappa$ where $\kappa$ is supercompact; and 
if we force further with $\Fn(\kappa,2)$, do we obtain a model of $\refl_{\RC}=\mu$ ?
\smallskip

[S.F., Ottenbreit Sakai III] should be mentioned here. 
\smallskip
}
\fi

It is consistent that $\continuum$ is ccc-g.\ supercompact, $\MA$ holds but a reflection 
principle with the reflection point $\LT\aleph_2$ weaker than \RC\ (i.e.\ the statement
$\refl_\RC=\aleph_2$) also holds:
Start from a model of \ZFC\ with two supercompact cardinals. Use the smaller supercompact 
to force Fodor-type Reflection Principle (\FRP, which is a reflection principle with the 
reflection point $\LT\aleph_2$, and \FRP\ follows from \RC). Then force by the standard ccc forcing for Martin's Axiom 
to make the larger supercompact cardinal (which survives the first extension) to make it 
the continuum. In \cite{fjetal}, it is shown that \FRP\ is preserved by ccc generic 
extension. Thus in the resulting model, we still have \FRP\ together with \MA\ and that the 
continuum is ccc-g.\ supercompact. 

\begin{Cor}{\rm (K\"onig \cite{koenig} see also \cite{sfetal-I})}
  Suppose that $\kappa$ is $\calP$-g.\ supercompact where $\calP$ is the class of all
  $\sigma$-closed \pos. Then any non-special tree $T$ has club 
  many non-special subtrees of size $<\kappa$. 
\end{Cor}
\prf $\sigma$-closed \pos\ preserve non-specialty of trees (Todor\v{c}evi\'c, see 
\cite{stevo1983}). \qedofCor

\begin{Cor}{\rm (Diagonal Reflection Principle, see \cite{cox}, \cite{sfetal-I})}\quad 
  \Label{P-gen-7} Let 
  \begin{xitemize}
  \item[] 
    $\calS:=\setof{\pairof{M,\seqof{S_a}{a\in M}}}{M\not=\emptyset,\,
      S_a\subseteq [M]^{\aleph_0}\mbox{ for all }a\in M}$.
  \end{xitemize}

  For
  $\pairof{M,\seqof{S_a}{a\in M}}$, $\pairof{N,\seqof{S_a}{a\in N}}\in\calS$, let
  \begin{xitemize}
  \item[] 
    $\pairof{M,\seqof{S^M_a}{a\in M}}\sqsubseteq_\calS\pairof{N,\seqof{S^N_a}{a\in N}}$\qquad $:\Leftrightarrow$
    \\[\jot]
    \mbox{}\hfill$M\subseteq N$, and $S^M_a=S^N_a\cap[M]^{\aleph_0}$ for all $a\in M$. 
  \end{xitemize}

  Let the property $P$ be defined by stipulating that $P$ holds 
  in $\pairof{M,\seqof{S_a}{a\in M}}\in\calS$ if and only if 
  $S_a$ is a stationary subset of $[M]^{\aleph_0}$ for all $a\in M$.

  Suppose that $\calP$ is a class of \pos\ \st\ all elements of $\calP$ are proper.  If  
  $\kappa$ is a $\calP$-g.\ supercompact, then we have $\refl_{stat}(\calS, P)\leq\kappa$. \qed
\end{Cor}

The inequality $\refl_{stat}(\calS, P)\leq\kappa$ in \Corabove\ is optimal in the following 
sense: Suppose that $\kappa$ is supercompact and $\mu<\kappa$ is \st\ there is a 
non-reflecting stationary set $S\subseteq[\mu]^{\aleph_0}$.\footnote{For this argument we 
  only need here the ``non-reflectingness'' of the sort that there are club many
  $X\in[\mu]^{\LT\mu}$ \st\ $S\cap[X]^{\aleph_0}$ is not  
  stationary in $[X]^{\aleph_0}$.} If $\poP:=\Fn(\kappa,2)$ and 
$\genG$ is a $(\uniV,\poP)$-generic filter, then $\kappa$ is $Cohen$-g.\ supercompact in
$\uniV[\genG]$. Since $S$ remains a non-reflecting stationary subset of $[\mu]^{\aleph_0}$ in $\uniV[\genG]$, we 
have $\uniV[\genG]\modelof{\mu<\refl_{stat}(\calS,P)\leq\kappa=\continuum}$.

On the other hand, it
is also consistent (modulo large cardinals) that
$\refl_{stat}(\calS,P)\linebreak<\kappa=\continuum$ holds for a ccc-g.\ supercompact cardinal $\kappa$: Suppose 
that $\kappa_0<\kappa$ are two  
supercompact cardinals and $\poP$ and $\genG$ are as above. Then, in $\uniV[\genG]$, 
$\kappa_0$ and $\kappa$ are both Cohen-g.\ supercompact. 
Thus, by \Corabove, we have
$\uniV[\genG]\models\refl_{stat}(\calS,P)\leq\kappa_0<\kappa=\continuum$.

\ifextended
{\extendedcolor\bigskip\mbox{}

The equality
$\uniV[\genG]\modelof{\mu<\refl_{stat}(\calS,P)}$ above can be seen using the following 
facts. All of them, possibly except \LemmaAof{P-gen-9}, \assertof{1}, are trivial and well-known. 

\begin{LemmaA}
  \Label{P-gen-8} Suppose $X\subseteq X'$ with $\cardof{X}\geq\aleph_1$. \wassertof{1}
  If $S$ is a stationary subset of $[X]^{\aleph_0}$, then 
  \begin{xitemize}
  \xitemA[x-gen-16-0] 
   $S':=\setof{x\in[X']^{\aleph_0}}{x\cap X\in S}$ 
  \end{xitemize}
  is a stationary subset of $[X']^{\aleph_0}$.\smallskip

 \wassert{2} If $S$ is a non stationary subset of $[X]^{\aleph_0}$ then $S'$ defined as 
 above is a non stationary subset of $[X']^{\aleph_0}$. 
\end{LemmaA}
\prf \assertof{1}: Suppose that $C'\subseteq[X']^{\aleph_0}$ is club. Let
\begin{xitemize}
\xitemA[x-gen-17] $C:=\setof{x\cap X}{\cardof{x\cap X}=\aleph_0,\ x\in C'}$. 
\end{xitemize}
\begin{Claim}
  \Label{Cl-gen-0}
  $C$ contains a club set $\subseteq[X]^{\aleph_0}$. 
\end{Claim}
\prfofClaim Let $\theta$ be sufficiently large and
$\calM=\pairof{\calH(\theta),\in, C', X', \sqsubseteq}$ where $C'$ and $X'$ are thought to be 
unary relations and $\sqsubseteq$ is a well-ordering on $\calH(\theta)$. Note that 
$\sqsubseteq$ introduces a build-in Skolem hull operator for the structure $\calM$. The 
Skolem hull operator will be denoted by $\sk_\calM(\cdot)$.

Let $C_0:=\setof{y\in[X]^{\aleph_0}}{\sk_\calM(y)\cap X=y}$. Then it is easy to see that
$C_0\subseteq C$ (note that $\sk_\calM(y)\cap X'\in C'$ for any $y\in[X]^{\aleph_0}$), and that 
$C_0$ is a club subset of $[X]^{\aleph_0}$. 
\qedofClaim\qedskip

By \Claimabove\ and since $S$ is stationary in $[X]^{\aleph_0}$, 
we have
$S\cap C\not=\emptyset$. Suppose $y\in S\cap C$. Then there is  
$x\in C'$ with $x\cap X=y$ by the definition \xitemof{x-gen-17} of $C$. We also have
$y\in S'$ by the definition \xitemof{x-gen-16-0} of $S'$. 
Thus $y\in S'\cap C'$. \smallskip

\assertof{2}: Suppose that $C$ is a club subset of $[X]^{\aleph_0}$ \st\
$C\cap S=\emptyset$. Then $C':=\setof{x\in[x']^{\aleph_0}}{x\cap X\in C}$ is a club subset 
of $[X']^{\aleph_0}$ and $C'\cap S'=\emptyset$. 
\qedofLemmaA
\begin{LemmaA}
  \Label{P-gen-9} \wassert{1} Suppose that there is a non-reflecting stationary  $E\subseteq E^\mu_{\omega}$ 
  subset of a regular $\mu>\aleph_1$. Then there is a stationary 
    $S\subseteq[\mu]^{\aleph_0}$ 
  \st\ 
  \begin{xitemize}
  \xitemA[x-gen-17-0] 
    $S\cap[X]^{\aleph_0}$ is non-stationary subset of $[X]^{\aleph_0}$ for any 
    $X\in[\mu]^{\LT\mu}$ with $\omega_1\subseteq X$.
  \end{xitemize}

  \wassert{2} For any $\mu>\aleph_1$, if $S\subseteq[\mu]^{\aleph_0}$ is a stationary 
  subset of $[\mu]^{\aleph_0}$ \st\ 
  $S\cap[X]^{\aleph_0}$ is non-stationary in $[X]^{\aleph_0}$ for all 
  $X\in[\mu]^{\LT\mu}$ with $\omega_1\subseteq X$, then, for any $\lambda\geq\mu$, 
  \begin{xitemize}
  \item[] 
    $S':=\setof{x\in[\lambda]^{\aleph_0}}{x\cap\mu\in S}$ 
  \end{xitemize}
  is a stationary subset of $[\lambda]^{\aleph_0}$ \st, for all 
  $X\in[\lambda]^{<\mu}$ with $\omega_1\subseteq X$, $S'\cap[X]^{\aleph_0}$ 
  is non-stationary in 
  $[X]^{\aleph_0}$. 
\end{LemmaA}

The proof of \LemmaAof{P-gen-9} uses the following Theorem by Shelah:
\begin{ThmA}
  \Label{P-gen-10} {\rm (Shelah \cite{sh794})} For any regular $\mu>\aleph_1$, if 
  $\seqof{E_\xi}{\xi<\mu}$ is a pairwise disjoint sequence of stationary subsets of
  $E^\mu_\omega$, then 
  \begin{xitemize}
  \xitemA[x-gen-18] 
    $S=\setof{a\in[\mu]^{\aleph_0}}{a\cap\omega_1\in\omega_1,\,\sup^+(a)\in E_{a\cap\omega_1}}$
  \end{xitemize}
  is a stationary subset of $[\mu]^{\aleph_0}$. Here, $\sup^+(a)$ for $a\subseteq\On$ 
  denotes the ordinal $\sup(\setof{\xi+1}{\xi\in a})$. \qed
\end{ThmA}

\noindent
\prfof{\LemmaAof{P-gen-9}} \assertof{1}: Let $E\subseteq E^\mu_\omega$ be a non-reflecting 
stationary subset of a regular $\mu>\aleph_1$. That $E$ is stationary in $\mu$ but 
$E\cap\alpha$ is not stationary in $\alpha$ for all $\alpha<\mu$ with $\cf(\alpha)>\omega$.

Let $\seqof{E_\xi}{\xi<\mu}$ be a partition of $E$ into stationary subsets of $\mu$ (such a 
partition exists always by a theorem of Solovay).

Let $S$ be the stationary subset of $[\mu]^{\aleph_0}$ defined by \xitemof{x-gen-18} for 
this sequence $\seqof{E_\xi}{\xi<\mu}$. $S$ is stationary by \Thmof{P-gen-10}. Thus, we are 
done by showing that this $S$ satisfies \xitemof{x-gen-17-0}.

Let $X\in[\mu]^{\LT\mu}$ be \st\ $\omega_1\subseteq X$. \smallskip

\noindent{\bf Case I.} $X$ has the maximal element, say $\alpha$. Note that all elements $a$ of $S$ 
do not have the maximal element (since otherwise $\sup^+(a)$ would be a successor ordinal
$\not\in E_{a\cap\omega_1}$). Thus, $\alpha$ is not covered by $S\cap[X]^{\aleph_0}$ and 
hence it is not stationary. \smallskip

\noindent{\bf Case II.} $\sup^+(X)$ is of cofinality $>\omega$.
Let $\alpha=\sup^+(X)$. By assumption there is a club $C\subseteq\alpha$ \st\
$C\cap E=\emptyset$. The set
\begin{xitemize}
\item[] $\setof{a\in[X]^{\aleph_0}}{\sup^+(a)\in C}$
\end{xitemize}
is a club in $[X]^{\aleph_0}$ disjoint from $S\cap[X]^{\aleph_0}$.
\smallskip

\noindent{\bf Case III.} $\sup^+(X)$ is of cofinality $\omega$. 
Let $\alpha=\sup^+(X)$. 

If $\alpha\not\in E$ then 
Then \begin{xitemize}
\item[] $\setof{a\in [X]^{\aleph_0}}{\sup^+(a)=\alpha}$
\end{xitemize}
is a club in $[X]^{\aleph_0}$ disjoint from $S\cap[X]^{\aleph_0}$.

Otherwise, there is $\xi<\mu$ \st\ $\alpha\in E_\xi$. In this case,
\begin{xitemize}
\item[] 
  $\setof{a\in[X]^{\aleph_0}}{\sup^+(a)=\alpha,\,a\cap\omega_1>\xi}$ 
\end{xitemize}
is a club in $[X]^{\aleph_0}$ disjoint from $S\cap[X]^{\aleph_0}$. 
\smallskip

\assertof{2}: $S'$ is stationary by \Lemmaof{P-gen-8},\,\assertof{1}.
For any $X\in[\lambda]^{\LT\mu}$ with $\omega_1\subseteq X$, 
$S\cap[X\cap\mu]^{\aleph_0}$ is non stationary in $[X\cap\mu]^{\aleph_0}$. 
By \Lemmaof{P-gen-8},\,\assertof{2}, it follows that 
\begin{xitemize}
\item[] 
  $S'\cap[X]^{\aleph_0}=\setof{x\in[\lambda]^{\aleph_0}}{x\cap\mu\in S,\, s\subseteq X}$\\[\jot]
  $\phantom{S'\cap[X]^{\aleph_0}}=\setof{x\in[X]^{\aleph_0}}{x\cap(X\cap \mu)\in S\cap [X\cap\mu]^{\aleph_0}}$
\end{xitemize}
is non-stationary. 
\qedofLemmaA

}\fi 

\end{document}